\documentclass[twoside,notitlepage,11pt]{article}

\usepackage{amssymb}
\usepackage[leqno]{amsmath}
\usepackage{amsfonts}
\usepackage{amsopn}
\usepackage{amstext}
\usepackage{amsthm}
\usepackage[all]{xy}

\providecommand{\cal}{\mathcal}
\renewcommand{\Bbb}{\mathbb}

\newenvironment{pf}{\begin{proof}}{\end{proof}}

\newcommand{\Cee}{{\cal{C}}}
\newcommand{\Dee}{{\cal{D}}}
\newcommand{\Ef}{{\cal{F}}}

\newcommand{\Err}{{\Bbb{R}}}

\newcommand{\lam}{{\lambda}}
\newcommand{\al}{\alpha}

\newcommand{\eps}{\varepsilon}
\renewcommand{\phi}{\varphi}
\renewcommand{\rho}{\varrho}

\newcommand{\rest}{\restriction}

\newcommand{\ntr}{n\in\omega}

\newcommand{\loe}{\leqslant}
\newcommand{\goe}{\geqslant}

\newcommand{\subs}{\subseteq}
\newcommand{\sups}{\supseteq}
\newcommand{\nnempty}{\ne\emptyset}

\newcommand{\ovr}{\overline}
\renewcommand{\iff}{\Longleftrightarrow}

\newcommand{\cl}{\operatorname{cl}}

\newcommand{\w}{\operatorname{w}}

\newcommand{\id}{\operatorname{id}}

\newcommand{\by}{/}

\newcommand{\setof}[2]{\{#1\colon #2\}}

\newcommand{\sett}[2]{\{#1\}_{#2}}
\newcommand{\sn}[1]{\{#1\}} 
\newcommand{\dn}[2]{\{#1,#2\}} 
\newcommand{\pair}[2]{\langle #1, #2 \rangle} 
\newcommand{\map}[3]{#1\colon #2 \to #3} 
\newcommand{\img}[2]{#1[#2]} 
\newcommand{\inv}[2]{{#1}^{-1}[#2]} 

\newcommand{\Sig}{\Sigma}

\newcommand{\I}{{\ensuremath{\mathbb I}}}
\providecommand{\nat}{\omega}

\newcommand{\im}{\operatorname{im}}

\newcommand{\norm}[1]{\|#1\|}

\newcommand{\weakstar}{\ensuremath{{weak}^*}}

\newcommand{\el}{{\mathbb{L}}}

\newcommand{\ess}{\operatorname{ess}}
\newcommand{\econst}{\operatorname{E}}

\newcommand{\omn}{{\omega_1}}

\newtheorem{tw}{Theorem}[section]
\newtheorem{wn}[tw]{Corollary}
\newtheorem{lm}[tw]{Lemma}
\newtheorem{prop}[tw]{Proposition}

\theoremstyle{definition}

\theoremstyle{remark}

\newcommand{\alom}{\al<\omega_1}

\title{Linearly ordered compacta and Banach spaces with a projectional resolution of the identity}

\author{Wies\l aw Kubi\'s\\
Instytut Matematyki,
Akademia \'Swi\c etokrzyska\\
ul. \'Swi\c etokrzyska 15, 25-406 Kielce, Poland\\
E-mail: \texttt{wkubis@pu.kielce.pl}
}

\begin{document}
\maketitle

\begin{abstract}
We construct a compact linearly ordered space $K_\omn$ of weight $\aleph_1$, such that the space $C(K_\omn)$ is not isomorphic to a Banach space with a projectional resolution of the identity, while on the other hand, $K_\omn$ is a continuous image of a Valdivia compact and every separable subspace of $C(K_\omn)$ is contained in a $1$-complemented separable subspace. This answers two questions due to O. Kalenda and V. Montesinos.

\vspace{2mm}


\end{abstract}


\section{Introduction}

A subspace $F$ of a Banach space $E$ is {\em complemented} in $E$ if there exists a bounded linear projection $\map PEE$ such that $F=PE$. More precisely, we say that $F$ is {\em $k$-complemented} if $\norm P\loe k$.
A Banach space $E$ has the {\em separable complementation property} if every separable subspace of $E$ is contained in a complemented separable one. Typical examples of such spaces are Banach spaces with a countably norming Markushevich basis, which are called  {\em Plichko spaces} (see Kalenda \cite{Kalenda}). In case of Banach spaces of density $\aleph_1$, the property of being Plichko is equivalent to the existence of a {\em bounded projectional resolution of the identity}, i.e. a transfinite sequence of projections on separable subspaces satisfying some continuity and compatibility conditions (the precise definition is given below). In particular, every Plichko space of density $\aleph_1$ is the union of a continuous chain of complemented separable subspaces. 

A question of Ond\v rej Kalenda \cite[Question 4.5.10]{Kalenda} asks whether every closed subspace of a Plichko space is again a Plichko space. 
We describe a compact connected linearly ordered space $K_\omn$ of weight $\aleph_1$ which is an order preserving image of a linearly ordered Valdivia compact constructed in \cite{K_classR} and whose space of continuous functions is not Plichko. This answers Kalenda's question in the negative.

During the $34^{th}$ Winter School on Abstract Analysis (Lhota nad Rohanovem, Czech Republic, 14--21 January 2006), Vicente Montesinos raised the question whether every Banach space with the separable complementation property is isomorphic to a space with a projectional resolution of the identity.
We show that every separable subspace of $C(K_\omn)$ is contained in a $1$-complemented separable subspace, answering the above question in the negative.

On the other hand, we show that a Banach space of density $\aleph_1$ has a projectional resolution of the identity, provided it can be represented as the union of a continuous increasing sequence of separable subspaces $\sett{E_\al}{\alom}$ such that each $E_\al$ is $1$-complemented in $E_{\al+1}$. We apply this result for proving that every $1$-complemented subspace of a $1$-Plichko space of density $\aleph_1$ is again a $1$-Plichko space. This gives a partial positive answer to a question of Kalenda \cite[Question 4.5.10]{Kalenda}.

\section{Preliminaries}

We use standard notation and symbols concerning topology and set theory. For example, $\w(X)$ denotes the weight of a topological space $X$.
A compact space $K$ is called {\em $\aleph_0$-monolithic} if every separable subspace of $K$ is second countable. 
Given a surjection $\map fXY$, the sets $f^{-1}(y)$, where $y\in Y$, are called the {\em fibers} of $f$ or {\em $f$-fibers}.
The letter $\omega$ denotes the set of natural numbers $\{0,1,\dots\}$.
We denote by $\omega_1$ the first uncountable ordinal and we shall write $\aleph_1$ instead of $\omega_1$, when having in mind its cardinality, not its order type.
We denote by $|A|$ the cardinality of the set $A$.
A set $C\subs\omega_1$ is {\em closed} if $\sup_{\ntr}\xi_n\in C$ whenever $\setof{\xi_n}{\ntr}\subs C$ is increasing; $C$ is {\em unbounded} if $\sup C=\omega_1$.
Given an ordinal $\lam$, a sequence of sets $\sett{A_\al}{\al<\lam}$ will be called {\em increasing} if $A_\al\subs A_\beta$ for every $\al<\beta<\lam$, and {\em continuous} if $A_\delta=\bigcup_{\xi<\delta}A_\xi$ for every limit ordinal $\delta<\lam$.

In this note we deal with Banach spaces of density $\loe\aleph_1$. Recall that a Banach space $E$ has the {\em separable complementation property} if every separable subspace of $E$ is contained in a complemented separable one.
Fix a Banach space $E$ of density $\aleph_1$. A {\em projectional resolution of the identity} (briefly: {\em PRI}) in $E$ is a sequence $\sett{P_\al}{\alom}$ of projections of $E$ onto separable subspaces, satisfying the following conditions:
\begin{enumerate}
	\item[(1)] $\norm{P_\al}=1$.
	\item[(2)] $\al\loe \beta\implies P_\al = P_\al P_\beta = P_\beta P_\al$.
	\item[(3)] $E=\bigcup_{\alom}P_\al E$ and $P_\delta E =\cl(\bigcup_{\xi<\delta}P_\xi E)$ for every limit ordinal $\delta<\omn$.
\end{enumerate}
Weakening condition (1) to $\sup_{\alom}\norm{P_\al}<+\infty$, we obtain the notion of a {\em bounded projectional resolution of the identity}.
For a survey on the use of PRI's in nonseparable Banach spaces and for a historical background we refer to Chapter 6 of Fabian's book \cite{Fabian}.
A Banach space of density $\aleph_1$ is a {\em 1-Plichko space} if it has a projectional resolution of the identity. This is different from (although equivalent to) the original definition: see Definition 4.2.1 and Theorem 4.2.5 in \cite{Kalenda}.
A space isomorphic to a 1-Plichko space is called a {\em Plichko space} or, more precisely, a {\em $k$-Plichko space}, where $k\goe1$ is the constant coming from the isomorphism to a $1$-Plichko space.
In fact, a $k$-Plichko space of density $\aleph_1$ can be characterized as a space having a bounded PRI $\sett{P_\al}{\alom}$ such that $k\goe\sup_{\alom}\norm{P_\al}$ (see the proof of Theorem 4.2.4(ii) in \cite{Kalenda}).
Of course, every Plichko space has the separable complementation property.

Recall that a compact space $K$ is called {\em Valdivia compact} (see \cite{Kalenda}) if there exists an embedding $\map jK{[0,1]^\kappa}$ such that $\inv j{\Sig(\kappa)}$ is dense in $K$, where 
$\Sig(\kappa)=\setof{x\in [0,1]^\kappa}{|\setof{\al}{x(\al)\ne0}|\loe\aleph_0}$. Compact spaces embeddable into $\Sig(\kappa)$ are called {\em Corson compacta}.
By the result of \cite{KM}, a space of weight $\aleph_1$ is Valdivia compact if and only if it can be represented as the limit of a continuous inverse sequence of metric compacta with all bonding mappings being retractions -- a property analogous to the existence of a PRI in a Banach space. Valdivia compacta are dual to $1$-Plichko spaces in the following sense: if $K$ is a Valdivia compact then $C(K)$ is $1$-Plichko and if $E$ is a $1$-Plichko space then the closed unit ball of $E^*$ endowed with the \weakstar~topology is Valdivia compact. See \cite[Chapter 5]{Kalenda} for details.

Fix a Banach space $E$ of density $\aleph_1$. A {\em skeleton} in $E$ is a chain $\Cee$ of closed separable subspaces of $E$ such that $\bigcup\Cee=E$ and $\cl(\bigcup_{\ntr}C_n)\in\Cee$, whenever $C_0\subs C_1\subs\dots$ is a sequence in $\Cee$.
Given a skeleton $\Cee$, one can always choose an increasing sequence $\sett{E_\al}{\alom}\subs\Cee$ so that $E=\bigcup_{\alom}E_\al$ and $E_\delta=\cl(\bigcup_{\al<\delta}E_\al)$ for every limit ordinal $\alom$. In particular, every $C\in \Cee$ is contained in some $E_\al$.
We shall consider skeletons indexed by $\omn$, assuming implicitly that the enumeration is increasing and continuous.

\begin{lm}\label{fasfaqwqwt}
Let $E$ be a Banach space of density $\aleph_1$ and assume $\Cee,\Dee$ are skeletons in $E$. Then $\Cee\cap \Dee$ is a skeleton in $E$. More precisely, if $\Cee=\sett{C_\al}{\alom}$ and $\Dee=\sett{D_\al}{\alom}$ then there exists a closed and unbounded set $\Gamma\subs\omn$ such that $C_\al=D_\al$ for $\al\in \Gamma$.
\end{lm}

\begin{pf}
Let $\Gamma=\setof{\alom}{C_\al=D_\al}$. It is clear that $\Gamma$ is closed in $\omn$. Fix $\xi<\omn$. Since $C_\xi$ is separable, we can find $\al_1$ such that $C_\xi\subs D_{\al_1}$. Similarly, we can find $\al_2>\al_1$ such that $D_{\al_1}\subs C_{\al_2}$. Continuining this way, we obtain a sequence $\al_1<\al_2<\dots$ such that $C_{\al_{2n-1}}\subs D_{\al_{2n}}\subs C_{\al_{2n+1}}$ for every $\ntr$. Let $\delta=\sup_{\ntr}\al_n$. Then $C_\delta=D_\delta$, i.e. $\delta\in \Gamma$. Thus $\Gamma$ is unbounded. \end{pf}

\begin{wn} Assume $\sett{E_\al}{\alom}$ is a skeleton in a Plichko space $E$ of density $\aleph_1$. Then there exists a closed cofinal set $C\subs\omn$ such that $\sett{E_\al}{\al\in C}$ consists of complemented subspaces of $E$.
\end{wn}

Note that the skeleton $\sett{E_\al}{\al\in C}$ from the above corollary consists of $k$-complemented subspaces, where $k$ is the constant coming from a renorming of $E$ to a $1$-Plichko space.

\section{Linearly ordered compacta}

We shall consider linearly ordered compact spaces endowed with the order topology. Given such a space $K$, we denote by $0_K$ and $1_K$ the minimal and the maximal element of $K$ respectively. As usual, we denote by $[a,b]$ and $(a,b)$ the closed and the open interval with end-points $a,b\in K$. Given two linearly ordered compacta $K,L$, a map $\map fKL$ will be called {\em increasing} if
$x\loe y\implies f(x)\loe f(y)$
holds for every $x,y\in K$. It is straight to see that every increasing surjection is continuous. In particular, an order isomorphism is a homeomorphism.
We denote by $\I$ the closed unit interval of the reals.

The following lemma belongs to the folklore. The argument given below can be found, for example, in \cite[Lemma 2.1]{HJNR} and \cite[Prop. 5.7]{K_classR}.

\begin{lm}\label{irqppw}
Let $K$ be a linearly ordered space and let $a,b\in K$ be such that $a<b$. Then there exists a continuous increasing function $\map fK\I$ such that $f(a)=0$ and $f(b)=1$.
\end{lm}

\begin{pf} Since $X$ is a normal space, by the Urysohn Lemma, we can find a continuous function $\map hX\I$ such that $h(a)=0$ and $h(b)=1$. Modifying $h$, without losing the continuity, we may assume that $f(x)=0$ for $x\loe a$ and $f(x)=1$ for $x\goe b$. Define
$f(x)=\sup\setof{h(t)}{t\loe x}$.
Then $f$ is increasing, $f(a)=0$ and $f(b)=1$. It is straight to check that $f$ is continuous. 
\end{pf}

\begin{prop}\label{jwetpjp}
Let $K$ be a linearly ordered compact. Then the set of all increasing functions is linearly dense in $C(K)$.
\end{prop}

\begin{pf}
Fix $f\in C(K)$ and $\eps>0$. Then $K=J_0\cup\dots \cup J_{k-1}$, where each $J_i$ is an open interval such that the oscillation of $f$ on $J_i$ is $<\eps$. Choose $0_K=a_0<a_1<\dots<a_n=1_K$ such that for every $i<n$ either $[a_i,a_{i+1}]$ is contained in some $J_j$ or else $|[a_i,a_{i+1}]|=2$. By Lemma \ref{irqppw}, for each $i<n$ there exists an increasing function $\map{h_i}{[a_i,a_{i+1}]}\I$ such that $h_i(a_i)=0$ and $h_i(a_{i+1})=1$. Define
$$g(t)=(1-h_i(t))f(a_i) + h_i(t)f(a_{i+1})\quad\text{ for }t\in [a_i,a_{i+1}].$$
Then $\map gK{\Err}$ is a piece-wise monotone continuous function such that $\norm{f-g}<\eps$. Finally, piece-wise monotone functions are linear combinations of increasing functions.
\end{pf}

Given a continuous surjection of compact spaces $\map fXY$, we shall say that $C(Y)$ is {\em identified with a subspace of $C(X)$ via} $f$, having in mind the space
$\setof{\phi f}{\phi\in C(Y)}$, which is linearly isometric to $C(Y)$. In other words, $\phi\in C(X)$ is regarded to be a member of $C(Y)$ if and only if $\phi$ is constant on the fibers of $f$.
The next statement is in fact a reformulation of \cite[Prop. 5.7]{K_classR}, which says that a linearly ordered compact is the inverse limit of a sequence of ``smaller" linearly ordered compacta.

\begin{prop}\label{fqaiopjapf}
Assume $K$ is a compact linearly ordered space of weight $\aleph_1$. Then there exist metrizable linearly ordered compacta $K_\al$, $\alom$ and increasing quotients $\map{q_\al}K{K_\al}$ such that $\sett{C(K_\al)}{\alom}$ is a skeleton in $C(K)$, where $C(K_\al)$ is identified with a subspace of $C(K)$ via $q_\al$. Moreover, for each $\al<\beta<\omega_1$ there exists a unique increasing quotient $\map{q^\beta_\al}{K_\beta}{K_\al}$ such that $q_\al = q^\beta_\al q_\beta$.
\end{prop}

\begin{pf} In view of Proposition \ref{jwetpjp}, we can choose a linearly dense set $\Ef\subs C(K)$ consisting of increasing functions and such that $|\Ef|=\aleph_1$.
Let $\sett{\Ef_\al}{\alom}$ be a continuous increasing sequence of countable sets such that $\Ef=\bigcup_{\alom}\Ef_\al$. Let $E_\al$ be the closed linear span of $\Ef_\al$. Clearly, $\sett{E_\al}{\alom}$ is a skeleton in $C(K)$. 

Fix $\alom$. Define the relation $\sim_\al$ on $K$ as follows:
$$x\sim_\al y \iff (\forall\;f\in\Ef_\al)\; f(x)=f(y).$$
It is easy to check that $\sim_\al$ is an equivalence relation on $K$ whose equivalence classes are closed and convex. Let $K_\al=K\by\sim_\al$. Then $K_\al$ is a second countable linearly ordered compact space and $\sim_\al$ induces an increasing quotient map $\map{q_\al}{K}{K_\al}$. If $\al<\beta$ then $\sim_\al\sups \sim_\beta$, therefore there exists a unique (necessarily increasing and continuous) map $q^\beta_\al$ satisfying $q_\al = q^\beta_\al q_\beta$.

Finally, observe that $f\in E_\al$ if and only if $f$ is constant on the equivalence classes of $\sim_\al$. This shows that $q_\al$ identifies $C(K_\al)$ with $E_\al$.
\end{pf}

Fix a linearly ordered compact $K$. We denote by $\el(K)$ the set of all $p\in K$ such that $|[x,p]|=2$ for some $x<p$ in $K$. Such a point $x$ will be denoted by $p^-$. A point $p\in K$ will be called {\em internal} if it is not isolated from either of its sides, i.e. intervals $[x,p]$, $[p,y]$ are infinite for every $x<p<y$.
A point $p\in K$ will be called {\em external} in $K$ if it is not internal.
In other words, $p\in K$ is external iff either $p\in \dn{0_K}{1_K}$ or $p\in \el(K)$ or $p=q^-$ for some $q\in \el(K)$.
Observe that if $K$ is connected then the only external points are $0_K$ and $1_K$. If $K$ is second countable then the set of all external points of $K$ is countable.

Now fix $f\in C(K)$. We say that $p\in K$ is {\em irrelevant} for $f$ if one of the following conditions holds:
\begin{enumerate}
	\item[(1)] $p=0_K$ and $f\rest[p,b]$ is constant for some $b>p$,
	\item[(2)] $p=1_K$ and $f\rest[a,p]$ is constant for some $a<p$,
	\item[(3)] $0_K<p<1_K$ and $f\rest [a,b]$ is constant for some $a<p<b$.
\end{enumerate}
We say that $p$ is {\em essential} for $f$ if $p$ is not irrelevant for $f$.
We denote by $\ess(f)$ the set of all essential points of $f$. Observe that $\ess(f)$ is closed in $K$ and $f$ is constant on every interval contained in $K\setminus \ess(f)$. In fact, it is not hard to see that $f$ is constant on every interval of the form $[a,b]$ where $a,b\in\ess(f)$, $(a,b)\cap \ess(f)=\emptyset$ and $(a,b)\nnempty$.

Now assume that $X$ is a closed subset of a linearly ordered compact $K$. Define
$$\econst(K,X) = \setof{f\in C(K)}{\ess(f)\subs X}.$$
In other words, $\econst(K,X)$ is
the set of all $f\in C(K)$ which are constant on every interval of $K$ which has one of the following form: $[0_K,0_X]$, $[1_X,1_K]$, $[p^-,p]$, where $p\in \el(X)\setminus \el(K)$.
Clearly, $\econst(K,X)$ is a closed linear subspace of $K$.

\begin{prop}\label{w4egiwsjegj}
Assume $K$ is a linearly ordered compact and $E\subs C(K)$ is separable. Then there exists a closed separable subspace $X$ of\/ $K$ such that $E\subs \econst(K,X)$.
\end{prop}

\begin{pf}
By Proposition \ref{jwetpjp}, there exists a countable set $F\subs C(K)$ consisting of increasing functions such that $E$ is contained in the closed linear span of $F$. Let 
$$X=\cl\Bigl(\bigcup_{f\in F}\ess(f)\Bigr).$$
Then $F\subs \econst(K,X)$, therefore also $E\subs \econst(K,X)$.
It remains to show that $\ess(f)$ is separable for every $f\in F$.

Fix an increasing function $f\in C(K)$ and let $Y=\ess(f)$, $Z=\img fY=\img fK$. Observe that $f\rest Y$ is an increasing two-to-one map onto a second countable linearly ordered space $Z\subs\Err$. It is well known that in this case $Y$ is separable. For completeness, we give the proof.

Fix a countable dense set $D\subs\img fY$ which contains all external points of $\img fY$. Then $D'=Y\cap \inv fD$ is countable. We claim that $D'$ is dense in $Y$. For fix a nonempty open interval $U\subs Y$. If $f\rest U$ is not constant then $f(x)<f(y)$ for some $x,y\in U$ and hence there exists $r\in D$ such that $f(x)<r\loe f(y)$. Thus $f^{-1}(r)\subs D'\cap [x,y]\subs D'\cap U$. Now assume that $f\rest U$ is constant. Then $|U|\loe2$, because $f$ is two-to-one. Let $U=\dn xy$, where $x\loe y$. We may assume that $Y$ is infinite and consequently either $y<1_Y$ or $0_Y<x$. Assume $y<1_Y$ (the other case is the same). Then $y=p^-$ for some $p\in\el(Y)$. If $f(y)<f(p)$ then $t=f(y)$ is isolated from the right. Otherwise, $x=y\in\el(Y)$ and $f(x^-)<f(x)$, because $f$ is two-to-one. In the latter case, $t=f(x)$ is isolated from the left. Thus $t\in D$ and hence $U\subs D'$.
\end{pf}

\begin{prop}\label{wejgpjwge4}
Let $K$ be a linearly ordered space which is a continuous image of a Valdivia compact. Then $K$ is $\aleph_0$-monolithic, i.e. every separable subspace of $K$ is second countable.
\end{prop}

\begin{pf}
Suppose $X\subs K$ is closed, separable and not second countable. Then $\el(X)$ is uncountable. Let $X_1$ be the quotient of $X$ obtained by replacing each interval of the form $[a^-,a]$, where $a\in \el(X)$, by a second countable interval $I_a$ such that $[a^-,a]$ has an increasing map onto $I_a$. For example, one can define $I_a=\I$ if $[a^-,a]$ is connected and $I_a=\dn01$ otherwise.
Then $X_1$ is a non-metrizable increasing image of a linearly ordered Valdivia compact. On the other hand, $X_1$ is first countable. By the result of Kalenda \cite{Ka99}, $X_1$ is Corson compact. Finally, Nakhmanson's theorem \cite{Na} (see also \cite[Thm. IV.10.1]{A}) says that $X_1$ is metrizable, a contradiction. 
\end{pf}

\section{Main lemmas}\label{klkweafa}

Let $K$ denote the double arrow space, i.e. the linearly ordered space of the form $(\I\times\sn0)\cup(0,1)\times\sn1$ endowed with the lexicographic order. Then $K$ is compact in the order topology and admits a natural two-one-one increasing quotient $\map qK\I$. Example 2 of Corson \cite{Corson61} shows that $C(\I)$ is not complemented in $C(K)$, when embedded via $q$. Corson's argument can be sketched as follows. Suppose $\map P{C(K)}{C(\I)}$ is a projection. Then $C(K)$ is isomorphic to $C(\I)\oplus E$, where $E=C(K)\by C(\I)$. On the other hand, it is easy to check that $E$ is isomorphic to $c_0(\I)$, therefore $C(K)$ is weakly Lindel\"of. On the other hand, $C(K)$ is not weakly Lindel\"of, because $K$ is a non-metrizable linearly ordered compact (by Nakhmanson's theorem \cite{Na}). Taking a separable space $F\subs C(\I)$ instead of $C(\I)$, one can repeat the above argument to show that 
$C(\I)$ is not contained in a separable complemented subspace of $C(K)$. 

Corson's argument uses essentially topological properties of nonseparable Banach spaces. Below we prove a more concrete result, which requires a direct argument. We shall apply it in the proof of our main result.

\begin{lm}\label{kluczowy}
Assume\/ $\map \theta KL$ is an increasing surjection of linearly ordered compacta such that the set
$$Q=\setof{x\in L}{x\text{ is internal in $L$ and }|\theta^{-1}(x)|>1}$$
is somewhere dense in $L$. Then $C(L)$ is not complemented in $C(K)$, when identified with the subspace of $C(K)$ via $\theta$.
\end{lm}

\begin{pf}
For each $x\in L$ define $x^-=\min\theta^{-1}(x)$ and $x^+=\max\theta^{-1}(x)$.
Then each fiber of $\theta$ is of the form $[x^-,x^+]$, where $x\in L$. Recall that $\theta$ identifies $C(L)$ with the set of all $f\in C(K)$ which are constant on every interval $[x^-,x^+]$, where $x\in L$. 
Suppose $\map P{C(K)}{C(K)}$ is a bounded linear projection onto $C(L)$.
Fix $N\in\nat$ such that
$$-1+N/3\goe\norm P.$$
Given $p\in Q$, choose an increasing function $\chi_p\in C(K)$ such that $\chi_p(t)=0$ for $t\loe p^-$ and $\chi_p(t)=1$ for $t\goe p^+$. Let $h_p=P\chi_p$. There exists a (unique) function $\ovr h_p\in C(L)$ such that $h_p=\ovr h_p \theta$. Define
$$Q^-=\setof{q\in Q}{\ovr h_q(q)<2/3}\quad\text{ and }\quad Q^+=\setof{q\in Q}{\ovr h_q(q)>1/3}.$$
Then at least one of the above sets is somewhere dense.
Further, define 
$$U^-_p=(\ovr h_p)^{-1}(-\infty,2/3)\quad\text{ and }\quad U^+_p=(\ovr h_p)^{-1}(1/3,+\infty).$$

Suppose that the set $Q^-$ is dense in the interval $(a,b)$. Choose $p_0<p_1<\dots<p_{N-1}$ in $Q^-\cap (a,b)$ so that $p_i\in U^-_{p_0}\cap \dots\cap U^-_{p_{i-1}}$ for every $i<N$. This is possible, because each $p_i$ is internal in $L$.
Choose $f\in C(K)$ such that $0\loe f\loe 1$ and
$$f\rest \theta^{-1}(p_i) = \chi_{p_i}\rest \theta^{-1}(p_i)\quad \text{ for }i<N$$
and $f$ is constant on $[p^-,p^+]$ for every $p\in L\setminus \{p_0,\dots,p_{N-1}\}$.
The function $f$ can be constructed as follows. For each $i<N-1$ choose a continuous function $\map{\phi_i}{[p_i,p_{i+1}]}\I$ such that $\phi_i(p_i)=1$ and $\phi_i(p_{i+1})=0$. Define
$$f(t) =
\begin{cases}
0 &\qquad \text{ if } t< p_0,\\
\chi_i(t) &\qquad  \text{ if } t\in [p_i^-,p_i^+],\; i<N,\\
\phi_i\theta(t) &\qquad  \text{ if } t\in [p_i^+, p_{i+1}^-],\; i<N-1,\\
1 &\qquad  \text{ if } t>p_{N-1}.
\end{cases}
$$
Let $g=f-\sum_{i<N}\chi_{p_i}$. Then $g$ is constant on each interval of the form $\theta^{-1}(p)$, where $p\in L$. Indeed, if $p\notin \{p_0,\dots,p_{N-1}\}$ then all the functions $f,\chi_{p_0},\dots,\chi_{p_{N-1}}$ are constant on $\theta^{-1}(p)$. If $t\in[p_j^-,p_j^+]$ then
$$g(t) = f(t) - \sum_{i\loe j}\chi_{p_i}(t) = f(t) - (j-1) - \chi_{p_j}(t) = j-1,$$
because $f(t)=\chi_{p_j}(t)$.
It follows that $g\in C(L)$, i.e. $Pg=g$. Hence
$$Pf = Pg + P\Bigl(\sum_{i<N}\chi_{p_i}\Bigr) = g + \sum_{i<N}h_{p_i}.$$
Now choose $t\in U_{p_0}^-\cap\dots\cap U_{p_{N-1}}^-$ such that $t>p_{N-1}$. Note that $|(P f)(t)|\loe \norm P$, because $0\loe f\loe 1$. On the other hand, $h_{p_i}(t)<2/3$ for $i<N$ and consequently
\begin{align*}
-\norm P\loe (Pf)(t) &= g(t) +\sum_{i<N}h_{p_i}(t)
= f(t) - \sum_{i<N}\chi_{p_i}(t) + \sum_{i<N}h_{p_i}(t)\\
&= f(t) - N + \sum_{i<N}h_{p_i}(t)
< 1 - \sum_{i<N}(2/3)\\
&= 1-N/3 \loe -\norm P,
\end{align*}
which is a contradiction.

In case where the set $Q^-$ is nowhere dense, we use the fact that $Q^+$ must be somewhere dense and we choose a decreasing sequence $p_0>p_1>\dots>p_{N-1}$ in $Q^+$ so that $p_i\in U^+_{p_0}\cap\dots \cap U^+_{p_{i-1}}$ for $i<N$. Taking $t\in U^+_{p_0}\cap\dots \cap U^+_{p_{N-1}}$ with $t<p_{N-1}$ and considering a similar function $f$, we obtain
\begin{align*}
\norm P\goe (Pf)(t) &= f(t) - \sum_{i<N}\chi_{p_i}(t) + \sum_{i<N}h_{p_i}(t)
= f(t) + \sum_{i<N}h_{p_i}(t)\\
&> -1 + \sum_{i<N}(1/3)
= -1+N/3 \goe \norm P,
\end{align*}
which again is a contradiction.
\end{pf}

Recall that, given a compact space $X$ and its closed subspace $Y$, a {\em regular extension operator} is a linear operator $\map T{C(Y)}{C(X)}$ such that $T$ is positive (i.e. $Tf\goe 0$ whenever $f\goe 0$), $T1=1$ and $(Tf)\rest Y = f$ for every $f\in C(Y)$. Observe that in this case $\norm T=1$. The operator $T$ provides an isometric embedding of $C(Y)$ into $C(X)$ such that the image is a $1$-complemented subspace.

\begin{lm}\label{erigosgjag} Assume $X$ is a closed subset of a linearly ordered compact $K$. Then there exists a regular extension operator $\map T{C(X)}{C(K)}$ such that $\econst(K,X)\subs TC(X)$.
\end{lm}

\begin{pf} For each $a\in \el(X)$ choose a continuous increasing function $\map{h_a}{[a^-,a]}\I$ such that $h_a(a^-)=0$ and $h_a(a)=1$. Define
$$(Tf)(p) =
\begin{cases}
f(p) & \text{ if }p\in X,\\
f(0_X) & \text{ if }p<0_X,\\
f(1_X) & \text{ if }p>1_X,\\
(1-h_a(p))f(a^-) + h_a(p)f(a) \quad&\text{ if }p\in(a^-,a)\text{ for some }a\in\el(X).
\end{cases}$$
It is straight to check that $Tf\in C(K)$ for every $f\in C(X)$. Further, $(Tf)\rest X=f$, $T1=1$ and $T$ is positive, therefore it is a regular extension operator.
Note that $Tf$ is constant both on $[0_K,0_X]$ and $[1_X,1_K]$.
Finally, if $f\in \econst(K,X)$ then $T(f\rest X)$ is constant on each interval of the form $[a^-,a]$ where $a\in\el(X)\setminus\el(K)$, therefore $f=T(f\rest X)$. This shows that $\econst(K,X)$ is contained in the range of $T$.
\end{pf}

The above lemma implies that $C(K)$ has the separable complementation property, whenever $K$ is a linearly ordered $\aleph_0$-monolithic compact space.

\section{The space $K_\omn$}\label{wjegipjap}

Let $\pair K\loe$ be a linearly ordered compact space. Define the following relation on $K$:
$$x\sim y \iff [x,y] \text{ is scattered}.$$
It is clear that $\sim$ is an equivalence relation and its equivalence classes are closed and convex, therefore $K\by\sim$ is a linearly ordered compact space, endowed with the quotient topology and with the quotient ordering (i.e. $[x]_\sim \loe [y]_\sim\iff x\loe y$). In case where $K$ is dense-in-itself, the $\sim$-equivalence classes are at most two-element sets. Let $\map qK{K\by \sim}$ denote the quotient map. We call $q$ the {\em connectification} of $K$. In fact, $K\by\sim$ is a connected space, because if $[x]_\sim<[y]_\sim$, then setting $a=\max [x]_\sim$ and $b=\min [y]_\sim$, we have that $a<b$ and $a\not\sim b$, therefore $(a,b)\nnempty$ and $[x]_\sim<[z]_\sim<[y]_\sim$ for any $z\in(a,b)$. 

In \cite{K_classR}, a linearly ordered Valdivia compact space $V_\omn$ has been constructed, which has an increasing map onto every linearly ordered Valdivia compact. The space $V_\omn$ is 0-dimensional, dense-in-itself and has weight $\aleph_1$. Moreover, every clopen interval of $V_\omn$ is order isomorphic to $V_\omn$. In particular, every nonempty open subset of $V_\omn$ contains both an increasing and a decreasing copy of $\omn$.

\begin{tw}\label{eihgowigeo}
Let $K_\omn=V_\omn\by\sim$, where $\map q{V_\omn}{K_\omn}$ is the connectification of $V_\omn$. Then
\begin{enumerate}
	\item[(a)] $K_\omn$ is a connected linearly ordered compact of weight $\aleph_1$.
	\item[(b)] $K_\omn$ is a two-to-one increasing image of a linearly ordered Valdivia compact.
	\item[(c)] Every separable subspace of $C(K_\omn)$ is contained in a separable $1$-complemented subspace.
	\item[(d)] $C(K_\omn)$ does not have a skeleton of complemented subspaces; in particular it is not a Plichko space.
\end{enumerate}
\end{tw}

\begin{pf} Clearly, $K_\omn$ satisfies (a) and (b). 
For the proof of (c), fix a separable space $E_0\subs C(K_\omn)$. By Proposition \ref{w4egiwsjegj}, there exists a closed separable subspace $X$ of $K_\omn$ such that $E_0\subs \econst(K_\omn,X)$. By Lemma \ref{erigosgjag}, $\econst(K_\omn,X)$ is contained in a $1$-complemented subspace of $C(K_\omn)$, isometric to $C(X)$. By Proposition \ref{wejgpjwge4}, $K_\omn$ is $\aleph_0$-monolithic, therefore $X$ is second countable. This shows (c).

Suppose now that $C(K_\omn)$ has a skeleton $\Ef$ consisting of complemented subspaces. By Proposition \ref{fqaiopjapf} and by Lemma \ref{fasfaqwqwt}, there exists an increasing surjection $\map h{K_\omn}L$ such that $L$ is metrizable and $C(L)\in\Ef$, when identified with a subspace of $C(K_\omn)$ via $h$. Then $L$ is order isomorphic to the unit interval $\I$, being a connected separable linearly ordered compact. The set
$$A=\setof{x\in L}{0_L<x<1_L\text{ and }|h^{-1}(x)|>1}$$
is dense in $L$, because every non-degenerate interval of $K_\omn$ contains a copy of $\omn$. By Lemma \ref{kluczowy}, $C(L)$ is not complemented in $C(K_\omn)$.
This shows (d) and completes the proof.
\end{pf}

\section{Constructing compatible projections}

In this section we show that a Banach space of density $\aleph_1$ is $1$-Plichko if (and only if) it has a skeleton consisting of $1$-complemented subspaces.

\begin{lm}\label{klocki} Assume $E$ is a Banach space which is the union of a continuous chain $\sett{E_\al}{\al<\kappa}$ of closed subspaces such that
for every $\al<\kappa$, $E_\al$ is $1$-complemented in $E_{\al+1}$.
Then there exist projections $\map {S_\al}EE$, $\al<\kappa$, such that $\norm{S_\al}=1$, $\im S_\al=E_\al$ and $S_\al S_\beta= S_\al = S_\beta S_\al$ whenever $\al\loe\beta<\kappa$.
\end{lm}

\begin{pf}
We construct inductively projections $\setof{S_\al^\beta}{\al\loe\beta\loe\kappa}$ with the following properties:
\begin{enumerate}
	\item[(a)] $\map{S^\beta_\al}{E_\beta}{E_\al}$ has norm $1$ for every $\al\loe\beta$,
	\item[(b)] $S^\beta_\al S^\gamma_\beta=S^\gamma_\al$ whenever $\al\loe\beta\loe\gamma$. 
\end{enumerate}
We start with $S^0_0=\id_{E_0}$. Fix $0<\delta\loe\kappa$ and assume $S^\beta_\al$ have been constructed for each $\al\loe\beta<\delta$ and they satisfy conditions (a), (b). There are two cases:

{\bf Case 1:} $\delta=\gamma+1$. Using the assumption, fix a projection $\map T{E_{\gamma+1}}{E_\gamma}$ with $\norm T=1$ and define $S^\delta_\al= S^\gamma_\al T$ for every $\al<\delta$. Clearly both (a) and (b) are satisfied.

{\bf Case 2:} $\delta$ is a limit ordinal. Let $G=\bigcup_{\al<\delta}E_\al$. Then $G$ is a dense linear subspace of $E_\delta$. Define $\map {h_\al}{G}{E_\al}$ by setting
$$h_\al(x)=S^{\beta}_\al x ,\qquad\text{ where }\qquad\beta=\min\setof{\xi\in[\al,\delta)}{x\in E_\xi}.$$
Note that if $\al<\beta$, $x\in E_\beta$ and $\gamma\in[\beta,\delta)$ then 
\begin{equation}
S^\gamma_\al x =S^\beta_\al S^\gamma_\beta x =S^\beta_\al x ,
\tag{$*$}\end{equation}
because of (b) and by the fact that $S^\gamma_\beta\rest E_\beta=\id_{E_\beta}$. Using ($*$), it is easy to see that $h_\al$ is a linear operator. Clearly $\norm{h_\al}=1$, thus it can be uniquely extended to a linear operator $\map{S^\delta_\al}{E_\delta}{E_\al}$. Finally, $\norm{S^\delta_\al}=\norm{h_\al}=1$ and $S^\delta_\al$ is a projection onto $E_\al$. Thus (a) holds.

It remains to show (b). Fix $\al<\beta<\delta$. By continuity, it suffices to check that $S^\delta_\al x=S^\beta_\al S^\delta_\beta x$ holds for every $x\in G$. Fix $x\in G$ and find $\gamma\in[\beta,\gamma)$ such that $x\in E_\gamma$. We have
$$S^\beta_\al S^\delta_\beta x = S^\beta_\al S^\gamma_\beta x = S^\gamma_\al x = S^\delta_\al x.$$
Thus both conditions (a) and (b) hold. It follows that the construction can be carried out. 

Finally, define $S_\al:=S^{\kappa}_\al$. Clearly, $S_\al$ is a projection of $E$ onto $E_\al$ and $\norm{S_\al}=1$. If $\al<\beta<\kappa$ then 
$$S_\al S_\beta=S^{\kappa}_\al S^{\kappa}_\beta = S^\beta_\al S^{\kappa}_\beta S^{\kappa}_\beta = S^\beta_\al S^{\kappa}_\beta = S^{\kappa}_\al = S_\al$$
and of course $S_\beta S_\al = S_\al$, because $E_\al\subs E_\beta$. This completes the proof.
\end{pf}

\begin{wn}\label{giojeweg} Assume $E$ is a Banach space with a skeleton $\sett{E_\al}{\alom}$ such that $E_\al$ is $1$-complemented in $E_{\al+1}$ for every $\alom$. Then $E$ is a $1$-Plichko space.
\end{wn}

The following application of Lemma \ref{klocki} provides a partial positive answer to a question of Kalenda \cite[Question 4.5.10]{Kalenda}.

\begin{tw}\label{wefeiogfewiqgfj}
Assume $E$ is a $1$-Plichko Banach space of density $\aleph_1$. Then every $1$-comple\-men\-ted subspace of $E$ is $1$-Plichko.
\end{tw}

\begin{pf}
Let $\setof{P_\al}{\al<\omn}$ be a PRI on $E$ and let $\map QEE$ be a projection with $\norm Q=1$. Let $F:=\im Q$ and assume $F$ is not separable. Define $E_\al=\im P_\al$, $F_\al=F\cap E_\al$ and $R_\al=QP_\al$. Note that $\im R_\al=\img Q{E_\al}$.
We claim that the set $S=\setof{\delta<\omn}{\img Q{E_\delta}=F_\delta}$ is closed and unbounded in $\omn$. Indeed, define
$$\phi(\al)=\min\setof{\beta<\omn}{\img Q{E_\al}\subs E_\beta}$$
and observe that $\phi$ is well defined, since each $\img Q{E_\al}$ is separable and $E=\bigcup_{\xi<\omn}E_\xi$. Now, if $\delta<\omn$ is such that $\phi(\al)<\delta$ whenever $\al<\delta$, then 
$$\img Q{E_\delta}= \img{Q}{\cl(\bigcup_{\al<\delta} E_\al)}\subs\cl\bigcup_{\al<\delta}\img Q{E_\al}\subs \bigcup_{\al<\delta}E_{\phi(\al)}\subs E_\delta,$$
therefore $\delta\in S$. Thus $S$ is unbounded in $\omn$. Without loss of generality we may assume that $S=\omn$. 

Now observe that $R_\delta\rest F$ is a projection of $F$ onto $F_\delta$ and therefore $F_\delta$ is $1$-complemented in $F$. Thus in particular $F=\bigcup_{\xi<\omn}F_\xi$, where $\sett{F_\xi}{\xi<\omn}$ is a chain of separable subspaces of $F$ satisfying conditions 1, 2 of Lemma \ref{klocki}. By this lemma we get a PRI on $F$, which shows that $F$ is $1$-Plichko.
\end{pf}

A slight modification of the above proof gives the following

\begin{wn}
Assume $F$ is a complemented subspace of a Banach space $E$ of density $\aleph_1$. If $E$ has a skeleton consisting of complemented subspaces then so does $F$.
\end{wn}

\section{Final remarks and questions}

By Theorem \ref{eihgowigeo}, the following two properties of a Banach space $E$ of density $\aleph_1$ turn out to be different:
\begin{enumerate}
	\item[$C_k$:] $E$ is the union of an increasing sequence of separable $k$-complemented subspaces,
	\item[$CC_k$:] $E$ has a skeleton of separable $k$-complemented subspaces.
\end{enumerate}
Property $C_k$ is equivalent to the fact that every separable subspace of $E$ is contained in a separable $k$-complemented subspace. Every $k$-Plichko space has property $CC_k$ and property $CC_1$ is equivalent to the fact that $E$ is $1$-Plichko (Corollary \ref{giojeweg}).
If $E$ satisfies $CC_k$ and $F$ is an $l$-complemented subspace of $E$ then $F$ satisfies $CC_{kl}$, by the arguments from the proof of Theorem \ref{wefeiogfewiqgfj}.
We do not know whether $CC_k$ implies $k$-Plichko, in case where $k>1$.
We also do not know whether a closed subspace of a Plichko space necessarily has the separable complementation property. Finally, we do not know whether a $1$-complemented subspace of a Banach space with property $C_1$ necessarily has the separable complementation property.

\end{document}